\newif\ifspringer
\springertrue
\springerfalse

\ifspringer
\documentclass[graybox]{svmult}
\else
\documentclass[a4paper,11pt]{article}
\usepackage{theorem}         
\fi

\usepackage{amsmath}
\usepackage{amsfonts}
\usepackage{amssymb}         
\usepackage{amsopn}
\usepackage{graphicx}        
\usepackage{mathrsfs}		
\usepackage{todonotes}
\usepackage{mathptmx}       
\usepackage{makeidx}         
\usepackage{graphicx,epsfig,color}
\usepackage{multicol}        
\usepackage{multirow}
\usepackage[bottom]{footmisc}
\usepackage{url}
\usepackage{amsfonts}
\usepackage{wrapfig}
\usepackage{amsmath}
\usepackage{algorithm}
\usepackage{algorithmic}
\usepackage{soul}

\newcommand{\transpose}[1]{{#1}^\top}

\ifspringer
\makeindex             
\else

\theoremstyle{changebreak}                

\fi

\begin{document}

\ifspringer
\title*{Unassigned distance geometry and the Buckminsterfullerene\thanks{This paper was partly supported by the ANR PRCI project “MultiBioStruct”.}}

\titlerunning{Distance geometry and the Buckminsterfullerene}
\author{Leo Liberti} 
\authorrunning{Liberti}
\institute{Leo Liberti\at LIX CNRS Ecole Polytechnique, Institut Polytechnique de Paris, 91128 Palaiseau, France \\ \email{liberti@lix.polytechnique.fr}}

%
\maketitle
\else

\thispagestyle{empty}

\begin{center}
{\LARGE Unassigned distance geometry and the Buckminsterfullerene}
\par \bigskip
{\sc Leo Liberti${}^1$}
\par \bigskip
{\small\it LIX CNRS, \'Ecole Polytechnique, Institut Polytechnique de Paris, F-91128 Palaiseau, France} \\ Email:\url{liberti@lix.polytechnique.fr}
\par \medskip August 2024
\end{center}
\par \bigskip

\fi


\ifspringer
\abstract{%
  \else
  \begin{abstract}
    \fi
The Buckminsterfullerene is an inorganic molecule consisting of 60 carbon atoms in the shape of a soccer ball. It was used in \cite{liga} to showcase algorithms that find the correct shape of a protein from limited data (length of inter-atomic distances) without any further chemical experiment: in that case, by means of a complicated constructive heuristic based on genetic algorithms. In this paper we show that we can reconstruct the Buckminsterfullerene structure by means of mathematical programming, standard solver software, and little else.
\ifspringer
}
\else
\end{abstract} 
\fi

\section{Introduction}
\label{s:intro}
The authors of \cite{liga} state that ``\textit{ab initio} structure solution of nanostructured materials [that cannot be solved using crystallographic methods] is feasible using diffraction data in combination with distance geometry methods''. They then proceed to find the structure of Buckminsterfullerene (a highly symmetric molecule consisting of 60 carbon atoms \cite{senn}) using a heuristic method, called LIGA, that employs distance values in order to construct small rigid polytopes which are then combined into larger structures. The combination is achieved by means of a genetic algorithm approach that involve tournaments for selection of the fittest. Some backtracking also occurs.

The input to LIGA (number of atoms and list of all inter-atomic distance values) is obtained from powder diffraction data: we note that the graph underlying the distance values is not a part of the input. The experiment in \cite[Fig.~2c-d]{liga} suggests that the algorithmic reconstruction process yields some errors; but that it suffices to increase the multiplicities of the distance values in the list to obtain the correct shape. 

In this paper we show, using Mathematical Programming (MP) methods, that a correct structure of the Buckminsterfullerene can be obtained from the original distance list, without having to add multiple copies of each distance value, and without having to create a new algorithm. Most of the work is carried out on one of two MP formulations involved in the reconstruction problem, as explained below. The only algorithmic contribution is given by a ``vanilla'' MultiStart (MS) approach.

Our interest in the Buckminsterfullerene only extends so far as it constitutes a well-known and interesting instance of the Unassigned Distance Geometry Problem (UDGP), defined below. Many people have managed to reconstruct the shape of the Buckminsterfullerene: most of them have used heuristic optimization methods to find ``good'' minima of an energy formula derived from the partial differential equations that define the atomic force field \cite{schoenLJ1,schoenLJ2,cheng-C60}. We follow a different path, namely that of \cite{liga}, which takes as input the values of the inter-atomic distances of the molecule. 

\subsection{Two problems in distance geometry}
\label{s:dg}
We now introduce the Distance Geometry Problem (DGP) and its ``unassigned'' variant, the UDGP. The input of the DGP is an integer $K$ and a simple edge-weighted undirected graph $G=(V,E,d)$ where $d:E\to\mathbb{R}_+$. The DGP asks whether there exists a \textit{realization} $x:V\to\mathbb{R}^K$ of $G$ in $\mathbb{R}^K$ satisfying
\begin{equation}
  \forall \{i,j\}\in E\quad \|x_i-x_j\|_2^2 = d_{ij}^2.
  \label{dgp}
\end{equation}
In other words, the DGP asks us to draw the graph $G$ in $K$ dimensions such that edges are drawn as segments having length equal to the weight \cite{dgbook}. The DGP is \textbf{NP}-hard \cite{saxe79}, it is in \textbf{NP} if $K=1$ but it is unknown whether it is in \textbf{NP} for $K>1$ \cite{dgpinnp}. It has applications in designing time synchronization protocols in wireless networks, in the localization of mobile sensors, in the determination of protein structure from Nuclear Magnetic Resonance (NMR) experiments, and in many other fields. See \cite{dgp-sirev,vetterli} for more information about the DGP.

The UDGP removes the adjacency information from the graph edges: its input is a pair of integers $K,n$ and a list $\delta$ of $m$ positive scalars. It asks whether there exists an assignment $\alpha$ mapping $[m]=\{1,\ldots,m\}$ to the set $\bar{E}=\{\{i,j\}\;|\;i,j\in[n]\}$, and a realization $x:[n]\to\mathbb{R}^K$ such that
\begin{equation}
  \forall \ell\le m \ \mbox{ with }\ \alpha(\ell)=\{i,j\} \quad \|x_i-x_j\|_2^2 = \delta_{\ell}^2.
  \label{udgp}
\end{equation}
Although the UDGP is not as well known as the DGP, in many applications it is the problem that is actually closest to the experiments leading to the use of the DGP. For example, the determination of protein structures from distances involves the solution of a DGP, but the actual raw input from the NMR experiments is closer to a list of values $\delta$ than to a graph $G$. In fact, $G$ is inferred from $\delta$ because the assignment $\alpha$ is computed first, together with a very imprecise realization \cite{wuthrich_nobel}, see Sect.~\ref{s:proteins} below. 

The UDGP has applications in the study of DNA \cite{skiena,dakic}, where it is known as the \textit{turnpike} or \textit{beltway} or \textit{partial digest} problem. It is the natural modelling framework for the study of nanostructures \cite{dg-4or} such as the Buckminsterfullerene. It also arises in the context of protein structure \cite{wuthrich_nobel}, but in this setting $|\delta|$ (the length of the sequence $\delta$) is much smaller than the number of atom pairs $N=n(n-1)/2$, which yields very imprecise assignment functions $\alpha$. Conversely, it is often the case that all (or almost all) pairwise distances can be measured by NMR interaction in molecules that are much smaller than proteins (e.g.~the Buckminsterfullerene) \cite{liga}. In turn, this yields much more precise assignments. The UDGP is \textbf{NP}-hard even if $m=N$ \cite{skiena}.

Both DGP and UDGP can be defined for any $K$, and there exist many applications for $K\in\{1,2,3\}$, as well as some for larger values of $K$. In this paper we fix $K=3$, since the Buckminsterfullerene is realized in 3D space.

\subsection{Goal of this paper}
\label{s:goal}
In this paper we solve the UDGP instance known as ``Buckminsterfullerene'' \cite{kroto} by means of MP, a formal language for specifying and solving optimization problems. We recall an existing Mixed-Integer Nonlinear Programming (MINLP) formulation of the UDGP, and propose a new Mixed-Integer Quadratically Constrained Programming (MIQCP) formulation. These can only be solved for tiny inputs much smaller than the Buckminsterfullerene. We therefore decompose the decision on the assignment $\alpha$ and the realization $x$: we construct an approximate Mixed-Integer Linear Programming (MILP) formulation of the UDGP which can be solved to larger sizes. We obtain an approximate assignment $\alpha'$ and an approximate realization, which we discard. We then use $\alpha'$ to define an edge-weighted graph $G_{\alpha'}$, and obtain a more precise realization by solving the DGP on $G_{\alpha'}$. Finally, we show the perfectly reconstructed shape of the Buckminsterfullerene from a distance list $\delta$ of size $N$.

\section{MINLP formulations for the UDGP}
\label{s:formulations}
The UDGP aims at finding an assignment $\alpha:[m]\to \bar{E}$ and a realization $x$ consistent with each other. We note that realizations are defined as functions from $V$ to $\mathbb{R}^K$ that satisfy the distance equations \eqref{dgp}. Therefore $x$ can be represented by an $n\times K$ matrix, where $n=|V|$: every row $x_i$ of $x$ is a vector in $\mathbb{R}^K$ for $i\in V$. We can encode the assignment $\alpha$ by means of binary decision variables $y_{ij\ell}$ such that $y_{ij\ell}=1$ iff $\alpha(\ell)=\{i,j\}$. The assignment $\alpha$ satisfies the following constraints: each distance value $\delta_\ell$ for $\ell\le m$ is assigned to exactly one edge $\{i,j\}$, and each pair $\{i,j\}$ is assigned to at most one $\ell$ (some pairs may remain unassigned).

\subsection{Exact MINLP formulations}
These considerations yield the following ``natural'' MP formulation for the UDGP, proposed in \cite{udgp}:
\begin{equation}
  \left.\begin{array}{rrcl}
    \min\limits_{x,y} & \sum\limits_{i<j\le n} \sum\limits_{\ell\le m} y_{ij\ell} (\|x_i-x_j\|_2^2 &-& \delta_\ell^2)^2 \\
    \forall \ell\le m & \sum\limits_{i<j\le n} y_{ij\ell} &=& 1 \\
    \forall i<j\le n & \sum\limits_{\ell\le m} y_{ij\ell} &\le& 1 \\
    & x\in\mathbb{R}^{nK} \quad\land\quad y\in\{0,1\}^{Nm}
  \end{array}\right\}
  \label{minlp}
\end{equation}
Since Eq.~\eqref{minlp} has nonlinear terms and integer variables, it is a MINLP formulation with a quintic polynomial as objective.

We propose a new MIQCP formulation. It decreases the polynomial degree of Eq.~\eqref{minlp} from five to two.
\begin{equation}
  {\small
  \left.\begin{array}{rrcl}
    \min\limits_{x,y,z} &\sum\limits_{i<j\le n}\sum\limits_{\ell\le m} (z^+_{ij\ell} &+& z^-_{ij\ell})\\
    \forall i<j\le n, \ell\le m &  -z^-_{ij\ell} - M(1\!-\!y_{ij\ell}) &\le& \|x_i-x_j\|_2^2 - \delta_\ell^2  \le  z^+_{ij\ell} + M(1\!-\!y_{ij\ell}) \\
    \forall \ell\le m & \sum\limits_{i<j\le n} y_{ij\ell} &=& 1 \\
    \forall i<j\le n & \sum\limits_{\ell\le m} y_{ij\ell} &\le& 1 \\
    & x\in\mathbb{R}^{nK} &\land &y\in\{0,1\}^{Nm} \ \land\  z\ge 0,
  \end{array}\right\}
  }
  \label{minlp2}
\end{equation}
where $M = (\sum_{\ell\le m} \delta_\ell)^2$ (see \cite[Prop.~2.2]{cordone}). ``Sandwiching'' the distance constraints allows constraint activity whenever $\alpha(\ell)=\{i,j\}$, and trivialize the constraint otherwise. The slack variables $z^-,z^+$ allow the solver to find an approximately feasible solution with imprecise distances. 

\subsection{An approximate MILP formulation}
The only nonlinear part of Eq.~\eqref{minlp2} is the Euclidean distance term $\|x_i-x_j\|_2$ in the middle member of the distance constraints. We have:
\[ \|x_i-x_j\|_2 = \|x_i\|_2+\|x_j\|_2-2\langle x_i,x_j\rangle = \langle x_i,x_i\rangle + \langle x_j,x_j\rangle - \langle x_i,x_j\rangle.\]
Now we linearize the nonlinear (inner product) terms: we create an $n\times n$ variable matrix $X$, and we replace each term $\langle x_i,x_j\rangle$ by $X_{ij}$ (which, by commutativity, yields a symmetric matrix $X$). For all $i<j$ and $\ell$, we obtain the linearized distance constraint
\[ -z_{ij\ell}^--M(1-y_{ij\ell}) \le X_{ii}+X_{jj}-2X_{ij} - \delta_\ell^2 \le z_{ij\ell}^++M(1-y_{ij\ell}). \]
So far, this is the usual first step in the construction of the Semidefinite Programming (SDP) relaxation for MP formulations with quadratic terms. However, instead of imposing the positive semidefinite (PSD) constraint $X\succeq 0$, we restrict our attention to a polyhedral inner approximation thereof, that of \textit{Diagonally Dominant} (DD) matrices. An $n\times n$ matrix is DD if it satisfies:
\begin{equation}
  \label{dd}
  \forall i\le n \quad X_{ii} \ge \sum\limits_{j\not=i} |X_{ij}|.
\end{equation}
By Gershgorin's Circle Theorem \cite{gershgorin}, every eigenvalue $\lambda_i$ of a matrix $A$ is contained in an interval $[A_{ii}-\sum_{j\not=i} |A_{ij}|\,,\, A_{ii}+\sum_{j\not=i} |A_{ij}|]$. By the definition of DD matrices (Eq.~\eqref{dd}), if $X$ is DD, then the lower interval bound $X_{ii}-\sum_{j\not=i} |X_{ij}|$ is non-negative, which implies that every eigenvalue of $X$ is non-negative, which means that $X$ is PSD. We note that the converse need not hold: many PSD matrices are not DD.

The DD cone is polyhedral: its extreme rays are given by the rank-one matrices $e_i\transpose{e}_i$ and $(e_i\pm e_j)\transpose{(e_i\pm e_j)}$ \cite{barker2}. We therefore proceed to linearize the absolute value term in Eq.~\eqref{dd} by means of an additional $n\times n$ matrix variable $T$, obtaining the linear reformulation \cite[Thm.~3.9]{ahmadimajumdar}:
\[  \forall i\le n \quad X_{ii} \ge \sum_{j\not=i} T_{ij} \qquad \land \qquad -T\le X\le T. \]
We remark that linear programming over the DD cone is also known as DD Programming (DDP).

We now put everything back together to yield the following MIDDP formulation:
\begin{equation}
  {\small
  \left.\begin{array}{rrcl}
    \min\limits_{X,T,y,z} &\sum\limits_{i<j\le n}\sum\limits_{\ell\le m} (z^+_{ij\ell} &+& z^-_{ij\ell})\\
    \forall i<j\le n, \ell\le m &  -z^-_{ij\ell}\!-\!M(1\!\!-\!\!y_{ij\ell}) &\le& X_{ii}\!+\!X_{jj}\!-\!2X_{ij}\!-\!d_\ell^2\!\le\!z^+_{ij\ell}\!+\!M(1\!\!-\!\!y_{ij\ell}) \\
    \forall \ell\le m & \sum\limits_{i<j\le n} y_{ij\ell} &=& 1 \\
    \forall i<j\le n & \sum\limits_{\ell\le m} y_{ij\ell} &\le& 1 \\
    \forall i\le n & X_{ii} &\ge& \sum\limits_{j\le n:j\neq i} T_{ij} \\
    & -T &\le& X \le T \\
    & X,T\in\mathbb{R}^{n\times n} &\land &y\in\{0,1\}^{Nm} \ \land\  z\ge 0.
  \end{array}\right\}
  }
  \label{middp}
\end{equation}
We note that Eq.~\eqref{middp} is a MILP, for which there exist technically viable solvers even for relatively large sizes \cite{cplex221}. We also remark that, although this reformulation may not immediately appear simple, it is based on a set of well-known building blocks.

Since not every PSD matrix is DD, a DDP is an inner approximation of the corresponding SDP, which means that, potentially, the DDP might be infeasible even if the SDP is not. This is not an issue for the MIDDP formulation in Eq.~\eqref{middp} since it is always feasible: it suffices to satisfy the DD constraints first, and then choose large enough slack variables $z^+,z^-$ to accommodate the error.

In practice, solving Eq.~\eqref{middp} yields an $n\times n$ matrix $X$. There are two ways to obtain an $n\times K$ (approximate) realization of $\delta$.
\begin{enumerate}
\item \textit{Principal Component Analysis} (PCA). Let $X=P\Lambda\transpose{P}$ be a spectral decomposition of $X$. The diagonal matrix $\Lambda$ is non-negative because $X$ is DD and therefore PSD, so that $\sqrt{\Lambda}$ is a real matrix. Let $\Lambda^{K}$ be the matrix obtained by $\Lambda$ by zeroing all of the diagonal entries aside from the $K$ largest. Then $\hat{x}=P\sqrt{\Lambda^K}$ has the property that $\bar{X}=\hat{x}\transpose{\hat{x}}=P\Lambda^K\transpose{P}$ is a best rank-$K$ approximation of $X$.
\item \textit{Barvinok's ``naive algorithm''} generalized to $K$ dimensions \cite{barvinok_orl}. Proceeding as for PCA, consider the factor $F=P\sqrt{\Lambda}$ of $X$, then sample (componentwise) an $n\times K$ matrix $Z\sim \mathsf{Normal}(0,1/\sqrt{K})$. The $n\times K$ matrix $\check{x}=FZ$ has good proximity properties with the constraints $\|x_i-x_j\|_2=\bar{y}_{ij\ell}\delta_\ell$ with arbitrarily high probability (here $\bar{y}$ are the assignment variables for $\alpha$ having value $1$ after solving the MIDDP: they simply select the correct $\ell$ to match to the $\{i,j\}$).
\end{enumerate}
Either of these two ways would allow us to retrieve an approximate realization $x$ from the DD matrix $X$, but experiments in this yielded very wrong realizations. 

\section{A matheuristic for the UDGP}
\textit{Matheuristics} are heuristic algorithms based on MP formulations \cite{matheuristics}. Our matheuristic decomposes the decision on the assignment $\alpha$ and the realization $x$.

\subsection{The way of proteins and their scale}
\label{s:proteins}
This decomposition idea is not new: the DGP application to protein conformation derives the distance values from NMR experiments, which in fact yield an ambiguous assignment. This assignment is disambiguated in two phases. The first rests on chemical considerations. The second phase is based on a metaheuristic (usually simulated annealing) which returns an assignment together with a poor quality realization, which is discarded. This assignment, however, is used to define a weighted graph that becomes the input of a DGP. The solution of this DGP provides the realization of interest, and the shape of the protein.

There are two technical differences with the Buckminsterfullerene: powder diffraction data does not yield any assignment (not even an ambiguous one), and the first disambiguation phase does not apply because of the atomic uniformity of the Buckminsterfullerene: it consists entirely of carbon atoms.

The main difference w.r.t.~methodologies for computing the shape of proteins and of inorganic molecules such as the Buckminsterfullerene, however, is one of scale. Proteins are much larger than the range of distances that can be measured using NMR, which is around 5.5{\AA}. As a consequence, only very few distances can be obtained. The protein graphs that provide the DGP input by means of which we compute realizations are very sparse.

The Buckminsterfullerene --- and similar molecules --- are much smaller than proteins. As a consequence, we can often obtain lists $\delta$ of distance values that have length comparable with the total number $N$ of index pairs. Measurements are never free from errors, so having ``complete'' lists of distances does not mean that the list is correct. The authors of \cite{liga} imply that a source of error is that $\delta$ may contain \textit{more} distances than needed: for each given distance value, its number of occurrences in $\delta$ may be overestimated. We do not simulate this type of error here, but we do look at how the structure reconstruction changes as $|\delta|$ grows to reach $N$.

\subsection{The algorithm}
For the above reasons, we propose a different methodology, based on MP, to implement the decomposition between $\alpha$ and $x$.

The matheuristic we propose is as follows:
\begin{enumerate}
\item compute an assignment $\alpha$ and an approximate realization $x'$ using the MIDDP (MILP) formulation Eq.~\eqref{middp};
\item discard the realization $x'$
\item construct the weighted graph $G$ from $\alpha$
\item define a DGP instance on $G$ and solve it, to obtain a realization $x^\ast$
\item return $x^\ast$.
\end{enumerate}

\subsubsection{Solving the MIDDP}
Even if MILP solvers are technically very advanced, the MIDDP formulation in Eq.~\eqref{middp} is very challenging to solve at the instance scale yielded by the Buckminsterfullerene. We therefore solve it ``locally'', which, for a MILP solver, means setting a resource limit (either number of Branch-and-Bound nodes, or time). The solution returned by the solver is an approximate assignment.

\subsubsection{Solving the DGP}
The DGP consists in finding a realization $x$ satisfying Eq.~\eqref{dgp}. Many different formulations of the DGP have been discussed in the literature \cite{dgp-sirev,dgds}. In this case we use the quartic formulation:
\begin{equation}
  \label{quartic}
  \min\limits_{x\in\mathbb{R}^{nK}} \sum\limits_{\{i,j\}\in E} (\|x_i-x_j\|_2^2 - d_{ij}^2)^2.
\end{equation}
It is unconstrained and nonconvex in the $x$ variables. As for the MIDDP, solving Eq.~\eqref{quartic} exactly is out of the technological reach of global optimization solvers, at least at the scale of the Buckminsterfullerene. We therefore use a simple MS heuristic: at each iteration we sample a random starting point, then call the local optimization solver IPOPT \cite{ipopt} to identify a close local minimum, choosing the best over the allowed iterations. 

\section{Finding the shape of the Buckminsterfullerene}
We implemented our formulation using AMPL \cite{ampl}. We used a smaller $M$ value in solving the MIDDP than the one mentioned in \cite[Prop.~2.2]{cordone}, since it is far from being tight. We set $M=\sum_{\ell}\delta_{\ell}^2$ in our computational experiments.

We solved the MIDDP formulation using CPLEX 22.1 \cite{cplex221}. The local solver deployed on Eq.~\eqref{quartic} is IPOPT 3.4.12. All experiments were carried out on an Apple M1 Max with 64GB RAM running MacOS Ventura 13.6.4. All software was natively compiled for the ARM64 architecture.

CPLEX was given a time limit of 360s (six minutes --- extending to one hour yielded no appreciable difference). The MS algorithm used to solve the DGP in Eq.~\eqref{quartic} was allowed to run for 10 iterations. IPOPT was given no time limit. 

\begin{figure}[!ht]
  \begin{center}
    \includegraphics[width=3.5cm]{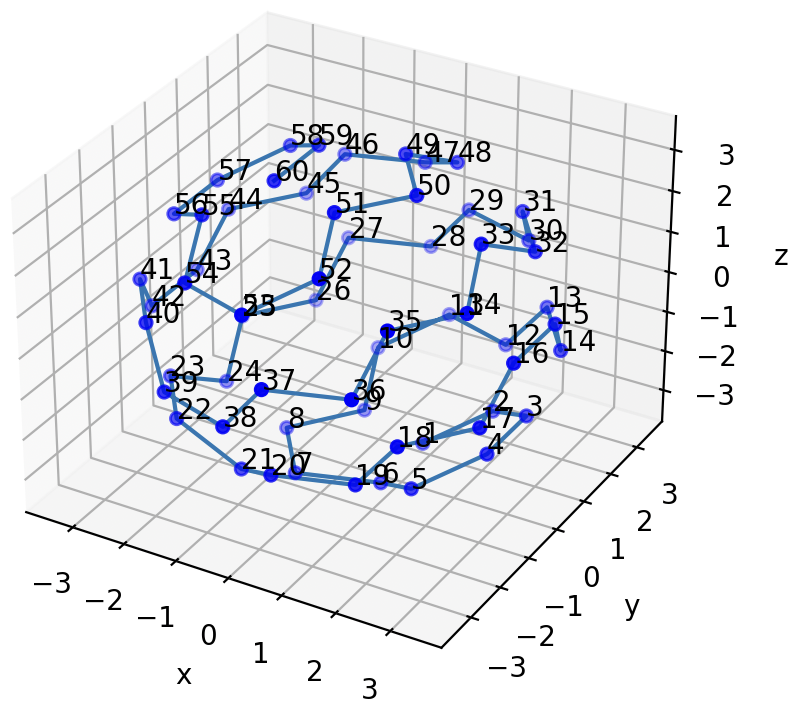}
    \includegraphics[width=3.5cm]{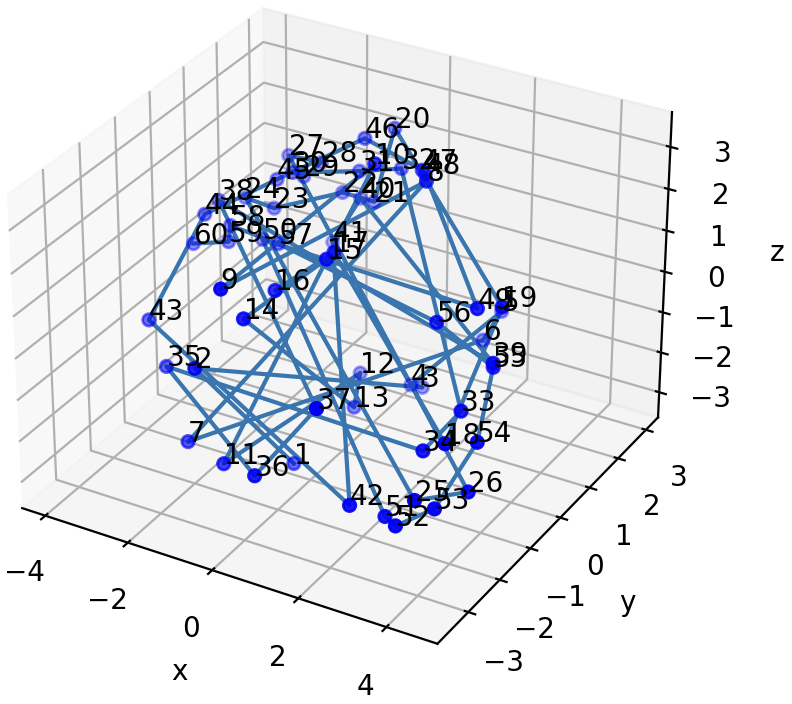}
    \includegraphics[width=3.5cm]{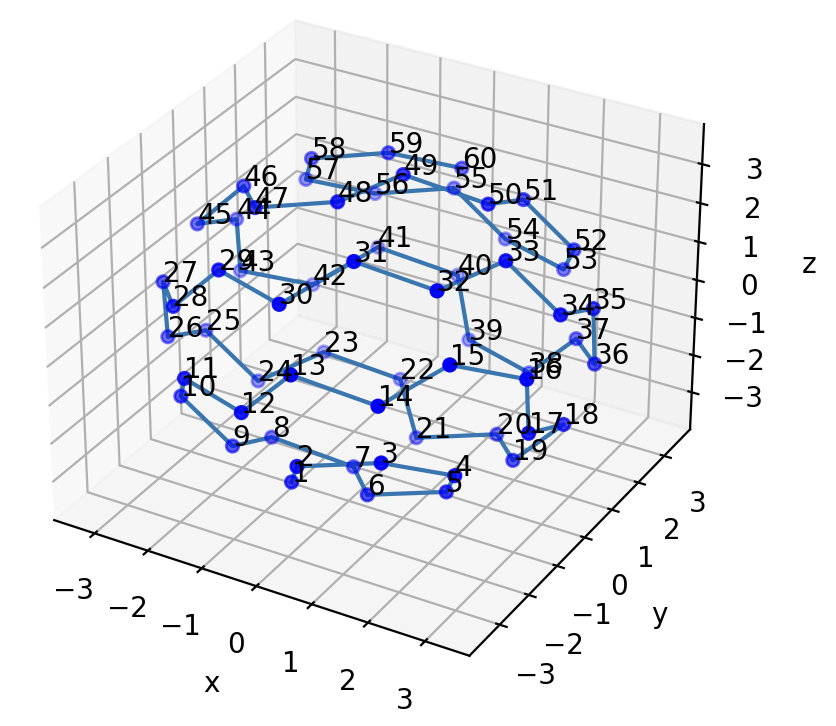}
  \end{center}
  \caption{Left: the Buckminsterfullerene molecule with a spiral order on the vertices. Center: wrong realizations obtained for $|\delta|=0.5\,N$. Right: the realization for $|\delta|=N$ yields a perfect reconstruction (rotated).}
  \label{f1}
\end{figure}

The Buckminsterfullerene has $n=60$ vertices and $N=1770$ distance values, many of them with multiplicities, given its molecular symmetry. We computed the exact realization according to \cite{senn} and extracted a complete list $\delta$ of its inter-atomic distances. From this instance we also generated instances with smaller $|\delta|$ through a probabilistic Erd\H{o}s-R\'enyi edge removal choice. We examined instances with $|\delta|\in[0.1\,N,1.0\,N]$, with $1.0\,N$ corresponding to the complete list. The original realization of the Buckminsterfullerene is given in Fig.~\ref{f1} (left). We also show a realization for $|\delta|=0.5\,N$ (Fig.~\ref{f1} center), which is rather typical of the whole range of probabilistically generated instances. Lastly, the perfectly reconstructed molecule corresponding to $|\delta|=1770$ is shown in Fig.~\ref{f1} (right). It may not appear identical because it is a rotation of the original: but the error returned by the solver was zero. 

\begin{figure}[!ht]
  \begin{center}
    \includegraphics[width=2.2cm]{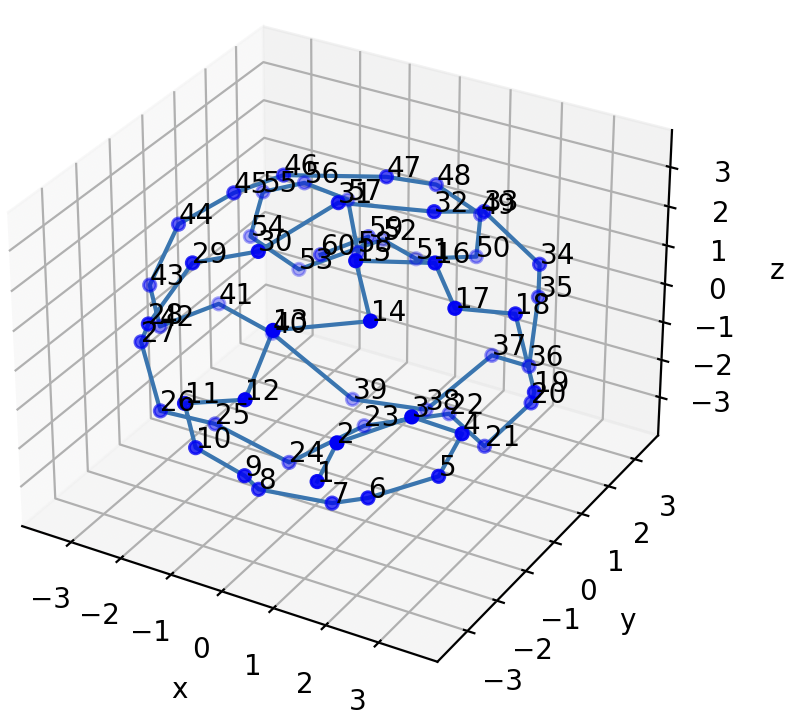}
    \includegraphics[width=2.2cm]{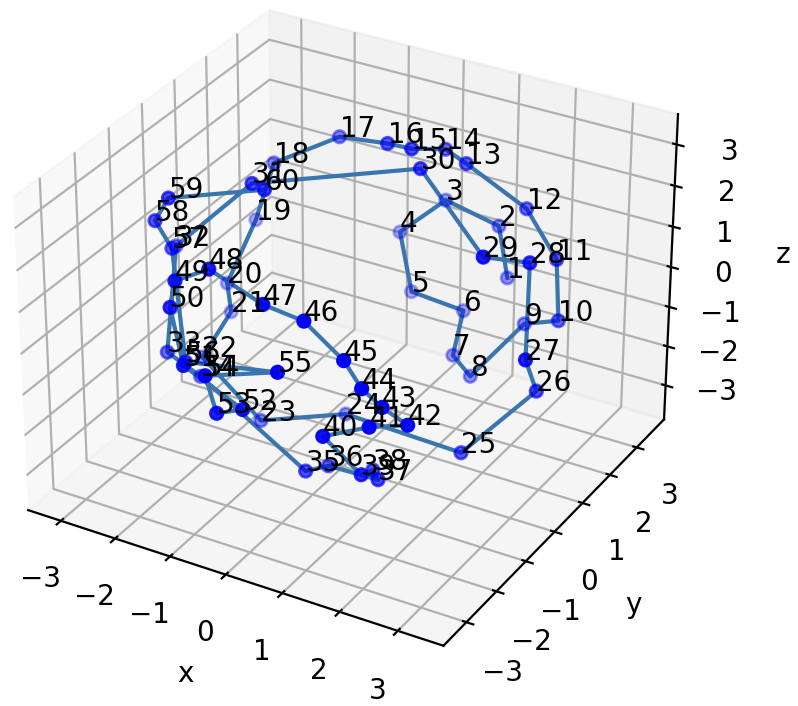}
    \includegraphics[width=2.2cm]{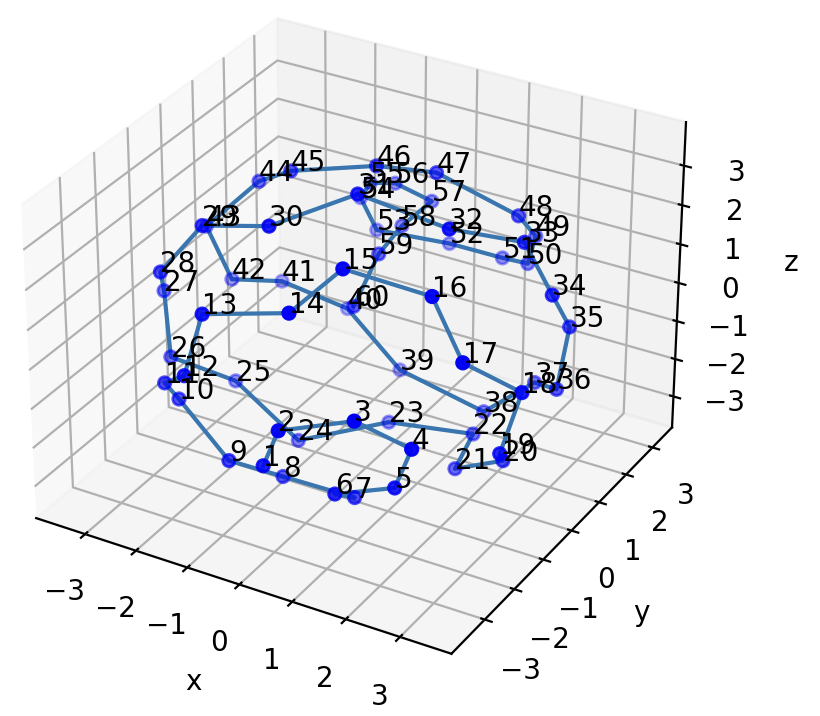}
    \includegraphics[width=2.2cm]{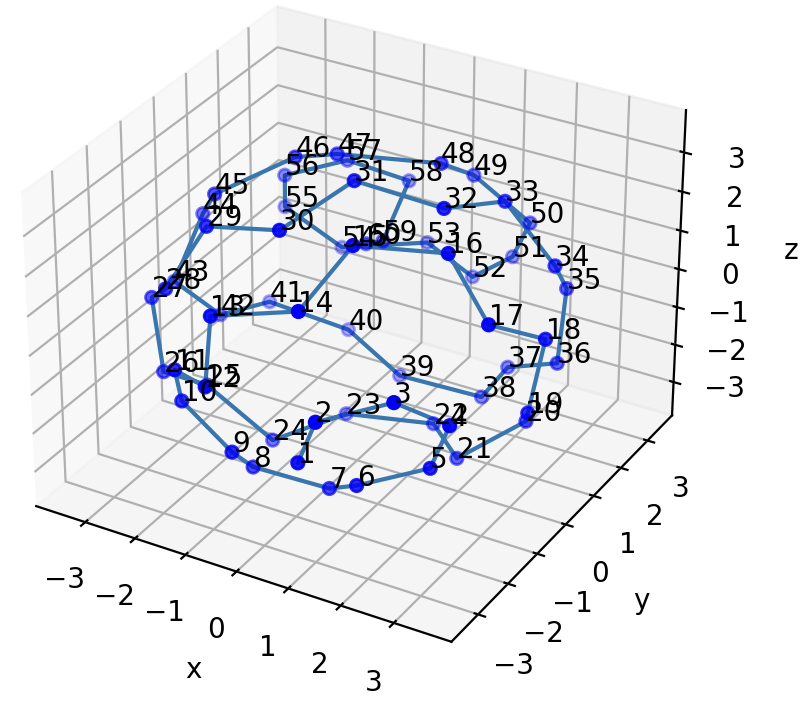}
    \includegraphics[width=2.2cm]{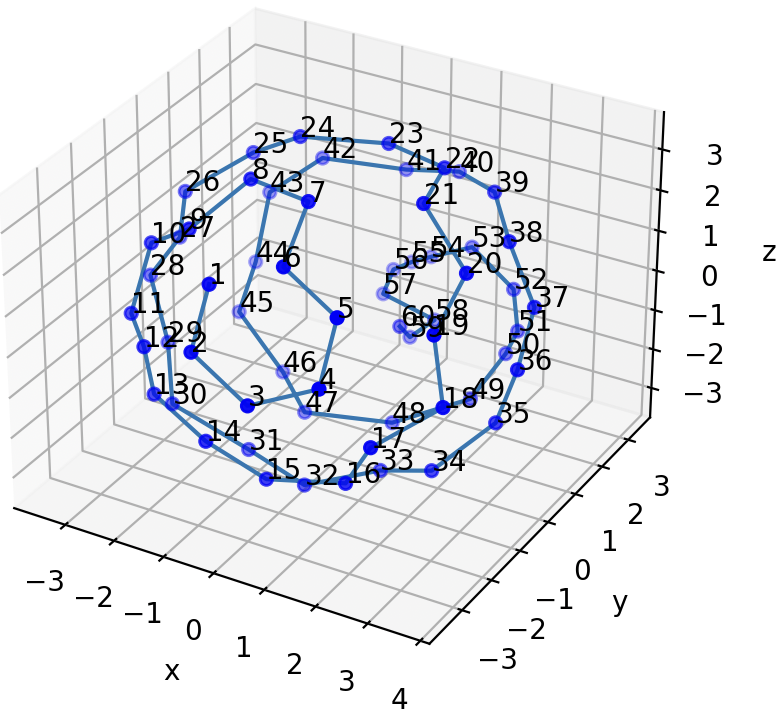}
  \end{center}
  \caption{Slightly wrong realizations for $|\delta|\in[0.995N,0.999N]$.}
  \label{f2}
\end{figure}

By trial-and-error we narrowed down an interval of interest for $|\delta|$, i.e.~when the resulting shape is not completely wrong but not yet correct, to $|\delta|\in[0.995\,N,0.999\,N]$. For this range we present numerical experiments based on the Mean Distance Error (MDE), i.e.~$(1/|E|)\sum_{\{i,j\}\in E} |\,\|x_i-x_j\|_2 - d_{ij}\,|$ computed on the reconstructed weighted graph.
  \begin{center}
    \begin{tabular}{l||r|r|r|r|r}
      density & 0.995 & 0.996 & 0.997 & 0.998 & 0.999 \\ \hline
          MDE & 0.722 & 1.057 & 0.778 & 0.552 & 0.603 
      \end{tabular}
  \end{center}
The realizations corresponding to above densities are shown in Fig.~\ref{f2}. As one can see, the realizations have the correct overall shape with some local errors.

We report the mean errors $(1/N)\|\mathsf{Adj}(G')-\mathsf{Adj}(G^\ast)\|_1$ of the reconstructed graphs $G'$ w.r.t.~the correct Buckminsterfullerene graph $G^\ast$, where $\mathsf{Adj}(G)$ is the upper triangular part of  weighted adjacency matrix of $G$.
\begin{center}
  \begin{tabular}{l||r|r|r|r|r}
    density & 0.995 & 0.996 & 0.997 & 0.998 & 0.999 \\ \hline
    mean error & 1.043 & 1.291 & 1.129 & 0.851 & 0.720 
  \end{tabular}
\end{center}

Finally, the CPU time taken for each instance to be solved was between 1000s and 1200s. 

\section{Conclusion}
The proposed methodology, simpler than that of \cite{liga}, allowed us to reconstruct the Buckminsterfullerene perfectly from the list of all its distances. It also allowed us to reconstruct locally wrong realizations with very few missing distances. This is in constrast to \cite{liga}, which required more than $N$ distances (same distance values but with higher multiplicities) in order to reconstruct the molecule perfectly.

\bibliographystyle{plain}
\bibliography{dr1}

\begin{thebibliography}{10}

\bibitem{ahmadimajumdar}
A.~Ahmadi and A.~Majumdar.
\newblock {DSOS} and {SDSOS} optimization: {M}ore tractable alternatives to sum
  of squares and semidefinite optimization.
\newblock {\em SIAM Journal on Applied Algebra and Geometry}, 3(2):193--230,
  2019.

\bibitem{barker2}
G.~Barker and D.~Carlson.
\newblock Cones of diagonally dominant matrices.
\newblock {\em Pacific Journal of Mathematics}, 57(1):15--32, 1975.

\bibitem{dgpinnp}
N.~Beeker, S.~Gaubert, C.~Glusa, and L.~Liberti.
\newblock Is the distance geometry problem in {NP}?
\newblock In A.~Mucherino, C.~Lavor, L.~Liberti, and N.~Maculan, editors, {\em
  Distance Geometry: Theory, Methods, and Applications}, pages 85--94.
  Springer, New York, 2013.

\bibitem{dg-4or}
S.~Billinge, P.~Duxbury, D.~Gon\c{c}alves, C.~Lavor, and A.~Mucherino.
\newblock Assigned and unassigned distance geometry: {A}pplications to
  biological molecules and nanostructures.
\newblock {\em 4OR}, 14:337--376, 2016.

\bibitem{cordone}
M.~Bruglieri, R.~Cordone, and L.~Liberti.
\newblock Maximum feasible subsystems of distance geometry constraints.
\newblock {\em Journal of Global Optimization}, 83:29--47, 2022.

\bibitem{cheng-C60}
Y.-H. Cheng, J.-H. Liao, Y.-J. Zhao, J.~Ni, and X.-B. Yang.
\newblock Theoretical investigations on stable structures of
  $\mathrm{C}_{60-n}\mathrm{N}_n$ ($n=2$--$12$): {S}ymmetry, model interaction,
  and global optimization.
\newblock {\em Carbon}, 154:140--149, 2019.

\bibitem{ipopt}
COIN-OR.
\newblock {\em Introduction to IPOPT: A tutorial for downloading, installing,
  and using IPOPT}, 2006.

\bibitem{dakic}
T.~Daki\'c.
\newblock {\em On the turnpike problem}.
\newblock PhD thesis, Simon Fraser University, 2000.

\bibitem{vetterli}
I.~Dokmani\'c, R.~Parhizkar, J.~Ranieri, and M.~Vetterli.
\newblock Euclidean distance matrices: Essential theory, algorithms and
  applications.
\newblock {\em IEEE Signal Processing Magazine}, 1053-5888:12--30, Nov. 2015.

\bibitem{schoenLJ2}
J.~Doye, R.~Leary, M.~Locatelli, and F.~Schoen.
\newblock Global optimization of {M}orse clusters by potential energy
  transformations.
\newblock {\em INFORMS Journal on Computing}, 16(4):371--379, 2004.

\bibitem{udgp}
P.~Duxbury, C.~Lavor, L.~Liberti, and L.~{de Salles-Neto}.
\newblock Unassigned distance geometry and molecular conformation problems.
\newblock {\em Journal of Global Optimization}, 83:73--82, 2022.

\bibitem{ampl}
R.~Fourer and D.~Gay.
\newblock {\em The {AMPL} Book}.
\newblock Duxbury Press, Pacific Grove, 2002.

\bibitem{gershgorin}
S.~Gershgorin.
\newblock \"uber die {A}bgrenzung der {E}igenwerte einer {M}atrix.
\newblock {\em Zvesti Akademii Nauk SSSR. Otdelenie Fizicheskikh i
  Matematicheskikh Nauk}, 6:749–754, 1931.

\bibitem{cplex221}
IBM.
\newblock {\em {ILOG CPLEX} 22.1 {U}ser's {M}anual}.
\newblock IBM, 2022.

\bibitem{liga}
P.~Juh\'as, D.~Cherba, P.~Duxbury, W.~Punch, and S.~Billinge.
\newblock \textit{Ab initio} determination of solid-state nanostructure.
\newblock {\em Nature}, 440(30):655--658, 2006.

\bibitem{kroto}
H.~Kroto, J.~Heath, S.~O'Brien, R.~Curl, and R.~Smalley.
\newblock $\mathrm{C}_60$: {B}uckminsterfullerene.
\newblock {\em Nature}, 318:162--163, 1985.

\bibitem{skiena}
P.~Lemke, S.~Skiena, and W.~Smith.
\newblock Reconstructing sets from interpoint distances.
\newblock In B.~Aronov and {\it et al.}, editors, {\em Discrete and
  Computational Geometry}, volume~25 of {\em Algorithms and Combinatorics},
  pages 597--631, Berlin, 2003. Springer.

\bibitem{dgds}
L.~Liberti.
\newblock Distance geometry and data science.
\newblock {\em TOP}, 28:271--339, 220.

\bibitem{dgbook}
L.~Liberti and C.~Lavor.
\newblock {\em Euclidean Distance Geometry: An Introduction}.
\newblock Springer, New York, 2017.

\bibitem{dgp-sirev}
L.~Liberti, C.~Lavor, N.~Maculan, and A.~Mucherino.
\newblock Euclidean distance geometry and applications.
\newblock {\em SIAM Review}, 56(1):3--69, 2014.

\bibitem{barvinok_orl}
L.~Liberti and K.~Vu.
\newblock Barvinok's naive algorithm in distance geometry.
\newblock {\em Operations Research Letters}, 46:476--481, 2018.

\bibitem{schoenLJ1}
M.~Locatelli and F.~Schoen.
\newblock Efficient algorithms for large scale global optimization:
  {L}ennard-{J}ones clusters.
\newblock {\em Computational Optimization and Applications}, 26:173--190, 2003.

\bibitem{matheuristics}
V.~Maniezzo, T.~St\"utzle, and S.~Vo{\ss}, editors.
\newblock {\em Matheuristics: {H}ybridizing metaheuristics and mathematical
  programming}, volume~10 of {\em Annals of Information Systems}, New York,
  2009. Springer.

\bibitem{saxe79}
J.~Saxe.
\newblock Embeddability of weighted graphs in $k$-space is strongly {NP}-hard.
\newblock {\em Proceedings of 17th Allerton Conference in Communications,
  Control and Computing}, pages 480--489, 1979.

\bibitem{senn}
P.~Senn.
\newblock Computation of the cartesian coordinates of {B}uckminsterfullerene.
\newblock {\em Journal of Chemical Education}, 72(4):302--303, 1995.

\bibitem{wuthrich_nobel}
K.~W\"uthrich.
\newblock {NMR} studies of structure and function of biological macromolecules
  ({N}obel lecture).
\newblock {\em Angewandte Chemie}, 42:3340--3363, 2003.

\end{thebibliography}
\end{document}